\documentclass[11pt]{article}

\usepackage{amsmath, epsfig, cite, lineno}
\usepackage{amssymb}
\usepackage{amsfonts, color}
\usepackage{latexsym}

\newtheorem{thm}{Theorem}

\newtheorem{lem}[thm]{Lemma}

\newcommand{\pf}{\noindent{\it Proof.} }

\def\inv{{\rm inv}}

\newcommand{\qed}{{\hfill\rule{4pt}{7pt}}\medskip}

\begin{document}


\begin{center}
{\Large\bf Bijective Proofs of Gould's and Rothe's Identities}
\end{center}

\vskip 2mm \centerline{Victor J. W. Guo\footnote{The author was
partially supported by a Junior Research Fellowship at the Erwin Schr\"odinger
International Institute for Mathematical Physics in Vienna.}}

\begin{center}
{\footnotesize Department of Mathematics, East China Normal
University, Shanghai 200062,
 People's Republic of China\\
{\tt jwguo@math.ecnu.edu.cn,\quad
http://math.ecnu.edu.cn/\textasciitilde{jwguo}}}
\end{center}


\vskip 0.7cm
\noindent{\bf Abstract.} We first give a bijective proof of
Gould's identity in the model of binary words. Then we
deduce Rothe's identity from Gould's identity again by a bijection,
which also leads to a double-sum extension of
the $q$-Chu-Vandermonde formula.


\section{Introduction}
There are two convolution formulas due to Rothe \cite{Rothe}:
\begin{align}
& \sum_{k=0}^{n}\frac{xy}{(x-kz)(y-(n-k)z)}{x-kz\choose k}{y-(n-k)z\choose n-k}=
\frac{x+y}{x+y-nz}{x+y-nz\choose n}, \label{eq:gou-1}\\[5pt]
& \sum_{k=0}^{n}\frac{x}{x-kz}{x-kz\choose k}{y+kz\choose n-k}=
{x+y\choose n}, \label{eq:gou-2}
\end{align}
which are famous in the literature. For example, Chu \cite{Chu}
used \eqref{eq:gou-1} and \eqref{eq:gou-2} to compute some determinants
involving binomial coefficients. For some generalizations of \eqref{eq:gou-1}
and \eqref{eq:gou-2}, we refer the reader to \cite{Schlosser,Strehl,Zeng}
and references therein.

Some proofs of \eqref{eq:gou-1} and \eqref{eq:gou-2}
can be found in \cite{Raney,Huang,Sprugnoli}.
It is not difficult to see that \eqref{eq:gou-1} can be deduced from \eqref{eq:gou-2}.
Blackwell and Dubins \cite{Blackwell} have given a
combinatorial proof of Rothe's identity \eqref{eq:gou-1},
which can also be proved in the model of lattice paths
(using \cite[p.~9]{Mohanty} or \cite[(1.1)]{Krattenthaler}).

Gould \cite{Gould,Gould57} reproved \eqref{eq:gou-1} and \eqref{eq:gou-2} and
also obtained the following interesting identity:
\begin{align}
\sum_{k=0}^{n}{x-kz\choose k}{y+kz\choose n-k}=
\sum_{k=0}^{n}{x+\epsilon-kz\choose k}{y-\epsilon+kz\choose n-k}. \label{eq:gou-3}
\end{align}

The main purpose of this paper is to give bijective proofs of Rothe's identity
\eqref{eq:gou-2} and Gould's identity \eqref{eq:gou-3} in the model of
binary words. As a conclusion, a double-sum extension of
the $q$-Chu-Vandermonde formula is also presented.

\section{Proof of \eqref{eq:gou-3}}
It is not difficult to see that Gould's identity \eqref{eq:gou-3} is equivalent to
\begin{equation}\label{eq:pqkm}
\sum_{k=0}^{n}{p-km\choose k}{q+km\choose n-k}
=\sum_{k=0}^{n}{p+1-km\choose k}{q-1+km\choose n-k}
\quad\text{for}\ p,q,m,n\in \mathbb{N}.
\end{equation}

We will prove that \eqref{eq:pqkm} holds for all integers $p\geq mn$
and $q\geq 1$ (and therefore for all real and complex numbers).

Let $\Gamma=\{a,b\}$ denote an alphabet with a grading $||a||=1$
and $||b||=m+1$. For a word $w=w_1\cdots w_n\in \Gamma^*$, its {\it length}
$n$ is denoted by $|w|$ and its {\it weight} by
$||w||=||w_1||+\cdots+||w_n||$. Let $|w|_b=(||w||-|w|)/m$ be the number of $b$'s
appearing in $w$, and let
$$
\Gamma_{p,k}:=\{w\in \Gamma^*\colon ||w||=p\ \text{and}\ |w|_b=k\}
\subseteq \Gamma^{p-km}.
$$
It is easy to see that $\#\Gamma_{p,k}={p-km\choose k}$. Furthermore, let
$$
\Gamma_{p,k}^{(r)}:=\{w\in \Gamma_{p,k}\colon
\text{$w$ has a prefix of weight $r$}\}.
$$

For $p,q\geq mn$, an obvious bijection (by factorization)
$$
\Gamma_{p+q,n}^{(p)}
\longleftrightarrow\biguplus_k \Gamma_{p,k}\times \Gamma_{q,n-k}
$$
leads to
$$
\#\Gamma_{p+q,n}^{(p)}
=\sum_{k}{p-km\choose k}{q-(n-k)m\choose n-k}.
$$
Thus, the identity \eqref{eq:pqkm} is equivalent to
\begin{equation}\label{eq:pq-num}
\#\Gamma_{p+q+mn,n}^{(p)}=\#\Gamma_{p+q+mn,n}^{(p+1)}.
\end{equation}

We need the following simple observation.
\begin{lem}\label{lem:ref}
Let $u,v\in \Gamma^*$ with $||u||,||v||\geq mn+1$, where
$n=|u\cdot v|_b$. Then there exist nonempty prefixes $x$
of $u$ and $y$ of $v$ such that $||x||=||y||$.
\end{lem}

\pf Suppose that $||u||=mn+r$ and $||v||=mn+s$ with $r,s\geq 1$.
Then the total number of nonempty prefixes of $u$ and $v$ is
$|u|+|v|=||u||+||v||-m|u\cdot v|_b=mn+r+s$.
On the other hand, each prefix of $u$ or $v$
has weight $\leq \max\{||u||,||v||\}$.
Hence $u$ and $v$ must have some nonempty prefixes of
the same weight.\qed

Now we can prove \eqref{eq:pq-num} by the following theorem.
\begin{thm}\label{thm:pq}
For $p\geq mn$ and $q\geq 1$, there is a bijection between
$\Gamma_{p+q+mn,n}^{(p)}$ and $\Gamma_{p+q+mn,n}^{(p+1)}$.
\end{thm}
\pf
Take any $w=u\cdot v\in \Gamma_{p+q+mn,n}^{(p)}$,
where $||u||=p$ and $||v||=q+mn$. Applying Lemma~\ref{lem:ref} to $v$ and
the reverse of $u\cdot a$, we see that $u$ has a suffix $x$ (possibly empty), i.e.,
$u=u'\cdot x$, and $v$ has a prefix $y$, i.e., $v=y\cdot v'$, such that
$||x||=||y||-1$.
Then $u'\cdot \overline y\cdot \overline x\cdot v'\in \Gamma_{p+q+mn,n}^{(p+1)}$,
where $\overline x$ and $\overline y$ are respectively the reverses of $x$ and $y$.
By selecting $x$ and $y$ with minimal length, we obtain a bijection.
\qed

\section{Proof of \eqref{eq:gou-2}}

Let us now consider $\Gamma_{p+q+mn,n}$ with $p\geq mn$ and $q\geq 1$. For each
$w\in \Gamma_{p+q+mn,n}$, let $w=u\cdot v$ denote the
unique factorization with $||u||\geq p$ but as small as possible.
There are two possibilities:
\begin{itemize}
\item If $||u||=p$, then $w\in \Gamma_{p+q+mn,n}^{(p)}$ and all
these words have been counted above.

\item If $||u||=p+j$ for some $1\leq j\leq m$, then the last letter
of $u$ must be a $b$. Namely, $u=u'\cdot b$ for some
$u'\in \Gamma_{p+j-m-1,k-1}$. The corresponding $v$
belongs to $\Gamma_{q+mn-j,n-k}$. It is easy to see that
the mapping $w\mapsto (u',v)$ may be inverted.
\end{itemize}

Hence there is
a bijection
\begin{align}
\Gamma_{p+q+mn,n}\longleftrightarrow
\Gamma_{p+q+mn,n}^{(p)}\biguplus_{j=1}^m \biguplus_{k=1}^{n}
\Gamma_{p+j-m-1,k-1}\times \Gamma_{q+mn-j,n-k}, \label{eq:bij}
\end{align}
which means that
\begin{equation}\label{eq:kmx}
\sum_{k=0}^{n}\left({p-km\choose k}{q+km \choose n-k}+
\sum_{j=1}^{m}{p-km+j-1\choose k-1}{q+km-j \choose n-k}\right)
={p+q\choose n}.
\end{equation}
However, by \eqref{eq:pqkm}, for $1\leq j\leq m$, we have
\begin{equation}\label{eq:kmpink}
\sum_{k=0}^{n}{p-km+j-1\choose k-1}{q+km-j \choose n-k}=
\sum_{k=0}^{n}{p-km-1\choose k-1}{q+km \choose n-k}.
\end{equation}
Combining \eqref{eq:kmx} and \eqref{eq:kmpink} yields
\begin{equation*}
\sum_{k=0}^{n}\left({p-km\choose k}+m{p-km-1\choose k-1}\right)
{q+km \choose n-k}={p+q\choose n},
\end{equation*}
or
\begin{equation*}
\sum_{k=0}^{n}\frac{p}{p-km}{p-km\choose k}
{q+km \choose n-k}={p+q\choose n},
\end{equation*}
which is Rothe's identity \eqref{eq:gou-2}.

\section{A new extension of the $q$-Chu-Vandermonde formula}
In this section we give a $q$-analogue of \eqref{eq:kmx}.
Recall that {\it $q$-binomial coefficient ${x\brack k}$}
is defined as $\prod_{i=1}^{k}(1-q^{x-i+1})/(1-q^i)$ if
$k\geq 0$ and $0$ otherwise. By \cite[Theorem 3.6]{Andrews}, we have
\begin{equation}\label{eq:invw}
\sum_{w\in \Gamma_{p,k}} q^{\inv(w)}=
{p-km\brack k}\quad (p\geq km),
\end{equation}
where $\inv (w)$ denotes the number of inversions of $w$.
Taking $\inv(w)$ into account and
using \eqref{eq:invw}, the bijection
\eqref{eq:bij} (replacing $p$ and $q$ by $x$ and $y$, respectively)
further implies that
\begin{align}
\sum_{k=0}^{n}q^{k(km+k+y-n)}\left({x-km\brack k}{y+km \brack n-k}+
\sum_{j=1}^{m}{x-km+j-1\brack k-1}{y+km-j \brack n-k}q^{-kj}\right)
={x+y\brack n}, \label{eq:q-chu}
\end{align}
which reduces to the $q$-Chu-Vandermonde formula if $m=0$, and reduces to
\begin{align*}
&\hskip -3mm
\sum_{k=0}^{n}q^{k(2k+y-n)}\left({x-k\brack k}{y+k \brack n-k}+
{x-k\brack k-1}{y+k-1\brack n-k}q^{-k}\right)={x+y\brack n}
\end{align*}
if $m=1$.

However, our bijection in Theorem~\ref{thm:pq} does not lead to the
corresponding $q$-analogue of \eqref{eq:kmpink}, and we cannot
simplify \eqref{eq:q-chu} to obtain a $q$-analogue of Rothe's
identity as before.

\vskip 2mm
\noindent{\bf Acknowledgments.}
The author is indebted to one of the referees for his
detailed constructive comments and suggestions, which enable the author
to shorten this note more than half. Especially, the present form of
Lemma~\ref{lem:ref} is his and it unifies two previous lemmas in this paper.

\end{document}